\documentclass{gtart_h}

\def\ifplaintex{\expandafter\ifx\csname documentclass\endcsname\relax}

\def\gtp{{\mathsurround=0pt\it $\cal G\mskip-2mu$eometry \&\ 
$\cal T\!\!$opology $\cal P\!$ublications}}  

\def\recd{{\small Received:\qua\receiveddate\ifx\reviseddate\relax
\else\qquad Revised:\qua\reviseddate\fi\par}} 

\def\lognumber#1{\def\thelognumber{#1}}
\def\volumenumber#1{\def\thevolumenumber{#1}}
\def\volumeyear#1{\def\thevolumeyear{#1}}
\def\papernumber#1{\def\thepapernumber{#1}}
\def\pagenumbers#1#2{\def\startpage{#1}\def\finishpage{#2}}
\def\published#1{\def\publishdate{#1}}

\def\received#1{\def\receiveddate{#1}}

\def\accepted#1{\def\accepteddate{#1}}

\def\asciiaddress#1{\def\theasciiaddress{#1}}
\def\asciiemail#1{\def\theasciiemail{#1}}

\long\def\asciiabstract#1{\long\def\theasciiabstract{#1}}

\let\\\par\let\thelognumber\relax\let\thevolumenumber\relax
\let\thepapernumber\relax\let\thevolumeyear\relax\let\startpage\relax
\let\finishpage\relax\let\publishdate\relax\let\receiveddate\relax
\let\reviseddate\relax\let\accepteddate\relax\let\theasciititle\relax
\let\theasciiauthors\relax\let\theasciiaddress\relax
\let\theasciiabstract\relax

\let\theasciiemail\relax

\ifplaintex
\font\logobig=cmssbx10 scaled 3836
\font\logomed=cmssbx10 scaled 2557
\else
\font\logobig=cmssbx10 scaled 4200
\font\logomed=cmssbx10 scaled 2800
\fi

\long\def\makeagttitle{   
\count0=\startpage
\agt\hfill      
\hbox to 45truept{\vbox to 0pt{\vglue -13truept{\logomed A\kern -.37em{\logobig 
T}\kern -.38em G}\vss}\hss}
\break
{\small Volume \thevolumenumber\ (\thevolumeyear)
\startpage--\finishpage\nl
Published: \publishdate}

\vglue .25truein

{\parskip=0pt\leftskip 0pt plus
1fil\def\\{\par\smallskip}{\Large\bf\thetitle}\par\medskip} \vglue
0.05truein

{\parskip=0pt\leftskip 0pt plus 1fil\def\\{\par}{\sc\theauthors}
\par\medskip}%
 
\vglue 0.03truein 

{\small\leftskip 25truept\rightskip 25truept{\bf Abstract}\stdspace\theabstract

{\bf AMS Classification}\stdspace\theprimaryclass
\ifx\thesecondaryclass\relax\else; \thesecondaryclass\fi\par
{\bf Keywords}\stdspace \thekeywords\par}\vglue 7truept

}   

\ifplaintex
\font\phead=cmsl9 scaled 950
\font\pnum=cmbx10 scaled 913
\font\pfoot=cmsl9 scaled 950
\headline{\vbox to 0pt{\vskip -4.5mm\line{\small\phead\ifnum
\count0=\startpage ISSN 1472-2739 (on-line) 1472-2747 (printed)
\hfill {\pnum\folio}\else\ifodd\count0\def\\{ }%
\ifx\theshorttitle\relax\thetitle\else\theshorttitle\fi\hfill{\pnum\folio}
\else\def\\{ and }{\pnum\folio}\hfill\ifx\theshortauthors\relax\theauthors
\else\theshortauthors\fi\fi\fi}\vss}}
\footline{\vbox to 0pt{\vglue 0mm\line{\small\pfoot\ifnum\count0=\startpage
\copyright\ \gtp\hfill\else
\agt, Volume \thevolumenumber\ (\thevolumeyear)\hfill\fi}\vss}}
\else
\font=cmsl9 scaled 1050
\font\lnum=cmbx10 
\font=cmsl9 scaled 1050
\makeatletter
\def\@oddhead{{\small\lhead\ifnum\count0=\startpage ISSN 1472-2739 
(on-line) 1472-2747 (printed)\hfill {\lnum\number\count0}\else\ifodd\count0
\def\\{ }\ifx\theshorttitle\relax \thetitle \else\theshorttitle\fi\hfill
{\lnum\number\count0}\else\def\\{ and }{\lnum\number\count0}
\hfill\ifx\theshortauthors\relax 
\theauthors\else\theshortauthors\fi\fi\fi}}\def\@evenhead{\@oddhead}
\def\@oddfoot{\small\lfoot\ifnum\count0=\startpage\copyright\ \gtp\hfill\else
\agt, Volume \thevolumenumber\ (\thevolumeyear)\hfill\fi}
\def\@evenfoot{\@oddfoot}
\makeatother
\fi
\let\maketitlepage\makeagttitle
\let\makeshorttitle\maketitlepage
\let\maketitle\maketitlepage

\newwrite\gtoutfile
\long\gdef\makeheadfile{  
{\def\\{, }\def\s{ }
\immediate\openout\gtoutfile head.xxx
\immediate\write\gtoutfile{Proxy-for: \ifx\theasciiauthors\relax
\theauthors\else\theasciiauthors\fi\s<\ifx\theasciiemail\relax\theemail\else\theasciiemail\fi>}
\immediate\write\gtoutfile{\noexpand\\}
\immediate\write\gtoutfile{Authors: \ifx\theasciiauthors\relax
\theauthors\else\theasciiauthors\fi}
{\def\\{ }\immediate\write\gtoutfile{Title: \ifx\theasciititle\relax
\thetitle\else\theasciititle\fi}}
\immediate\write\gtoutfile{Subj-class: GT or SG, GR etc}
\immediate\write\gtoutfile{MSC-class: \theprimaryclass\ifx\thesecondaryclass\relax\else, \thesecondaryclass\fi}
\immediate\write\gtoutfile{Journal-ref: Algebr. Geom. Topol. \thevolumenumber\s
(\thevolumeyear) \startpage-\finishpage}
\immediate\write\gtoutfile{Comments: Published by Algebraic and
Geometric Topology at}
\immediate\write\gtoutfile{\s\s\s  http://www.maths.warwick.ac.uk/agt/AGTVol\thevolumenumber/agt-\thevolumenumber-\thepapernumber.abs.html}
\immediate\write\gtoutfile{\noexpand\\}
\immediate\write\gtoutfile{}
\ifx\theasciiabstract\relax
\immediate\write\gtoutfile{\theabstract}\else
\immediate\write\gtoutfile{\theasciiabstract}\fi
\immediate\write\gtoutfile{}
\immediate\write\gtoutfile{\noexpand\\}
\immediate\write\gtoutfile{}
\immediate\closeout\gtoutfile}}  

\def\maketitlepage{\makeagttitle\makeheadfile}
\let\makeshorttitle\maketitlepage
\let\maketitle\maketitlepage

\lognumber{69}
\volumenumber{5}
\volumeyear{2005}
\papernumber{69}
\pagenumbers{1711}{1718}
\received{1 June 2005} 
\accepted{7 November 2005}
\published{7 December 2005}
\usepackage{amssymb,amsmath,amscd} 
\let\Bbb\mathbb

\newtheorem{theorem}{Theorem}[section]
\newtheorem{corollary}[theorem]{Corollary}

\newtheorem{proposition}[theorem]{Proposition}
\newtheorem{lemma}[theorem]{Lemma}

\theoremstyle{definition}

\newtheorem{question}[theorem]{Question}

\newtheorem{problem}[theorem]{Problem}

\newcommand{\edim}{{\rm edim}}
\newcommand{\z}{{\Bbb Z}}

\newcommand{\re}{{\Bbb R}}
\newcommand{\rp}{{\Bbb R}P}

\newcommand{\s}{{\Bbb S}}

\newcommand{\lo}{\longrightarrow}
\newcommand{\sm}{\setminus}

\begin{document}
\title{Extensions of maps to the projective plane}                    
\authors{Jerzy Dydak\\Michael Levin}                  
\address{Department of Mathematics, University of Tennessee\\
 Knoxville, TN 37996-1300, USA\\{\rm and}\\Department of Mathematics,
Ben Gurion University of the Negev\\
P.O.B. 653, Be'er Sheva 84105, ISRAEL}

\asciiaddress{Department of Mathematics, University of Tennessee\\
 Knoxville, TN 37996-1300, USA
\\and\\
Department of Mathematics,
Ben Gurion University of the Negev\\
P.O.B. 653, Be'er Sheva 84105, ISRAEL}

\gtemail{\mailto{dydak@math.utk.edu}{\rm\qua
and\qua}\mailto{mlevine@math.bgu.ac.il}}
\asciiemail{dydak@math.utk.edu, mlevine@math.bgu.ac.il}                     
\urladdr{http://www.math.utk.edu/~dydak}  

\begin{abstract}
It is proved that for a $3$-dimensional compact  metrizable space
 $X$ the infinite real projective space $\rp^\infty$ is an absolute extensor of $X$
if and only if the real projective plane $\rp^2$ is an absolute extensor of $X$
(see Theorems 1.2 and 1.5).
\end{abstract}

\asciiabstract{%
It is proved that for a 3-dimensional compact metrizable space X
the infinite real projective space is an absolute
extensor of X if and only if the real projective plane is an
absolute extensor of X.}

\maketitle

\begin{section}{Introduction}
Let $X$ be a compactum (= separable metric space)
and let $K$ be a CW complex. $K\in AE(X)$ (read: $K$ is an {\it absolute extensor} of $X$) or $X\tau K$
means that every map $f : A \lo K$, $A$ closed in $X$, extends over $X$.
  The {\it extensional dimension} $\edim X$ of $X$ is said
 to be dominated by a CW-complex $K$, written $\edim X \leq K$,
  if $X\tau K$. Thus for
  the covering dimension $\dim X$ of $X$ the condition
  $\dim X \leq n$ is equivalent to $\edim X \leq S^n$
  where $S^n$ is an $n$-dimensional sphere and for the cohomological
  dimension $\dim_G X$ of $X$ with respect to an abelian group
  $G$, the condition $\dim_G X\leq n$ is equivalent to
  $\edim X \leq K(G,n)$ where $K(G,n)$ is an Eilenberg-Mac Lane
  complex of type $(G,n)$.

Every time the coefficient group in homology is not explicitely stated, we
  mean it to be integers.

In case of CW complexes $K$ one can often reduce the relation $\edim X\leq K$
to $\edim X\leq K^{(n)}$, where $K^{(n)}$ is the $n$-skeleton of $K$.
 \begin{proposition}
  \label{skeleta}
Suppose $X$ is a compactum and $K$ is a CW complex.
If $\dim X\leq n$, then $\edim X\leq K$ is equivalent to $\edim X\leq K^{(n)}$.
 \end{proposition}
The proof follows easily using $\edim X\leq n$ to push maps off higher cells.

Dranishnikov \cite{dr2} proved the following
  important  theorems connecting  extensional and cohomological
  dimensions.

  \begin{theorem}
  \label{t-1}
  Let $K$  be a CW-complex and let a compactum $X$ be such that
  $\edim X \leq K$. Then
  $\dim_{H_n(K)} X \leq n$ for every $n>0$.
  \end{theorem}

   \begin{theorem}
  \label{t-2}
  Let $K$  be a simply connected CW-complex and let a compactum
  $X$ be finite
  dimensional.
  If $\dim_{H_n(K)} X \leq n$ for every $n>0$, then
  $\edim X \leq K$.
 \end{theorem}
 The requirement in Theorem \ref{t-2} that $X$ is finite dimensional
  cannot be omitted. To show that  take
 the famous infinite-dimensional compactum $X$ of
  Dranishnikov
   with $\dim_\z X  =3$ as in \cite{dr0}.
   Then the conclusion of Theorem \ref{t-2}
   does not hold for $K=S^3$.
  Let us mention in this connection another
  result \cite{l0}: there is a compactum $X$ satisfying the following
conditions:
\begin{itemize}
\item[(a)]  $\edim X > K$
   for every
   finite CW-complex $K$ with $\tilde H_*(K) \neq 0$,
\item[(b)]
   $\dim_G X \leq 2$ for every abelian group $G$,
\item[(c)] $\dim_G X
  \leq 1$    for every finite abelian group $G$.
\end{itemize}
   Here $\edim X > K$ means that
   $\edim X \leq K$ is false.

    With no restriction on $K$, Theorem \ref{t-2} does not hold.
  Indeed, the conclusion of Theorem \ref{t-2} is not satisfied
  if $K$ is a non-contractible  acyclic CW-complex and $X$ is
  the $2$-dimensional  disk.   Cencelj and  Dranishnikov \cite{cdr2}
  generalized  Theorem
  \ref{t-2}  for  nilpotent CW-complexes $K$ (see \cite{cdr1}
for the case of $K$ with fundamental group being finitely generated).

The real projective plane
  $\rp^2$ is the simplest CW-complex not covered by
   Cencelj-Dranishnikov's result. Thus we arrive at the following
   well-known open problem in Extension Theory.

   \begin{problem}
   \label{prob}
   Let $X$ be a finite dimensional compactum.
   Does $\dim_{\z_2} X \leq 1$ imply $\edim X \leq \rp^2$?
   More generally, does $\dim_{\z_p} X \leq 1$ imply
    $\edim X \leq M(\z_p,1)$, where $M(\z_p,1)$ is a Moore complex of
   type $(\z_p, 1)$?
   \end{problem}
   It is not difficult to  see that
    this problem can be  answered affirmatively
   if $\dim X \leq 2$ (use \ref{skeleta}). Sharing a belief that
   Problem \ref{prob} has a negative answer in higher
   dimensions the authors made
   a few unsuccessful  attempts to construct a counterexample
   in the first non-trivial case  $\dim X=3$
   and were  surprised to discover the following result.

 \begin{theorem}
 \label  {t1} Let $X$ be a compactum of dimension at most three.
If $\dim_{\z_2}X$ $\leq 1$, then $\edim X \leq \rp^2$.
 \end{theorem}

Notice (see \cite{dr3}) that there exist compacta $X$ of dimension 3
such that $\dim_{\z_2}X$ $\leq 1$, so Theorem \ref{t1} is not vacuous.

 This paper is devoted to proving of Theorem \ref{t1}. Theorem
 \ref{t1}
 can be formulated in a slightly different form.
 Let $X$ be a compactum.
 Take $\rp^\infty$
 as an Eilenberg-Mac Lane complex $K(\z_2,1)$.
  Then
  $\dim X \leq 3$   and  $\dim_{\z_2} X
 \leq 1$ imply $\edim X \leq \rp^3$. On the other hand
  by  Theorem \ref{t-1} the condition $\edim X
 \leq \rp^3$ implies that $\dim_{\z_2} X = \dim_{H_1(\rp^3)} X
 \leq 1$,  $\dim_\z X = \dim_{H_3(\rp^3)} X \leq 3 $ and
 by Alexandroff's theorem
 $\dim X \leq  3$ if $X$ is finite dimensional
(note that $H_k(\rp^n)=\z_2$ if $1\leq k < n$ is odd,
 $H_k(\rp^n)=0$ if $1 < k \leq n$ is even,
and $H_n(\rp^n)=\z$ if $n$ is odd - see p.89 of \cite{Wh}). Thus Theorem
 \ref{t1} is equivalent to the following, more general result.

 \begin{theorem}
 \label{t2} Let $X$ be a compactum of finite dimension.
 If $\edim X \leq \rp^3$, then $\edim X \leq \rp^2$.
  \end{theorem}
  We end this section with two questions related to
  Theorems \ref{t1} and \ref{t2}.

  \begin{question}
  \label{q1} Let $X$ be a compactum of dimension at most three.
 Does $\dim_{\z_p}X \leq 1$ imply  $\edim X \leq M(\z_p,1)$?
 \end{question}

\begin{question}
 \label{q2} Does $\edim X \leq \rp^3$ imply $\edim X \leq \rp^2$
for all, perhaps infinite-dimensional, compacta $X$?
  \end{question}
 \end{section}

 \begin{section}{Preliminaries}

\sh{Maps on projective  spaces}

Recall that the real projective $n$-space $\rp^n$ is obtained
from the $n$-sphere $S^n$ by identifying points $x$ and $-x$.
The resulting map $p_n:S^n\lo \rp^n$ is a covering projection
and $\rp^1$ is homeomorphic to $S^1$.
By $q_n:B^n\lo \rp^n$ we denote the quotient map
of the unit $n$-ball $B^n$ obtained by identifying $B^n$
with the upper hemisphere of $S^n$.
We consider all spheres to be subsets of the infinite-dimensional
sphere $S^\infty$. Similarly, we consider all projective
spaces $\rp^n$ to be subsets of the infinite projective space $\rp^\infty$.
Clearly, there is a universal covering projection $p:S^\infty\lo \rp^\infty$.
It is known that $\rp^\infty$ has a structure of a CW complex
making it an Eilenberg-MacLane complex of type $K(\z_2,1)$ as $S^\infty$ is contractible.

 \begin{proposition}
 \label  {p0} Any map $f : \rp^1 \lo
 \rp^2$ extends to a map $f' : \rp^2 \lo \rp^2$.
 \end{proposition}
 \proof It is obvious if $f$ is
 null-homotopic. Assume that $f$ is not homotopic to
 a constant map. Since $\pi_1(\rp^2)=\z_2$ and $\rp^1$ generates
 $\pi_1(\rp^2)$, $f$ is homotopic to the inclusion map of $\rp^1$ to $\rp^2$. Obviously,
that inclusion extends to the
 identity map of $\rp^2$, so $f$ extends over $\rp^2$.
 \endproof

  \begin{proposition}
 \label  {p1} If $f : \rp^2 \lo
 \rp^2$ induces the zero homomorphism of the fundamental groups, 
then $f$ extends to a map  $f': \rp^3 \lo \rp^2$.
 \end{proposition}
 \proof 
 Since $f$ induces the zero homomorphism of the
 fundamental group, $f$ can be lifted to $\beta : \rp^2 \lo S^2$.
 Since $H_2(\rp^2)=0$, the map $\gamma=\beta \circ q_3|_{\partial B^3} :
 \partial B^3  \lo S^2$
 induces the zero homomorphism $\gamma_* : H_2(\partial B^3) \lo H_2(S^2)$
  and hence $\gamma$  is null-homotopic. Thus $\gamma$ can be extended over $B^3$
  and this extension induces the corresponding extension of $f$ over
  $\rp^3$.
\endproof

 \begin{proposition}
 \label  {ExtendOverTorus} Let $Y$ be a topological space.
A map $f:S^1\times \rp^1\lo Y$ extends over $S^1\times \rp^2$
if and only if the composition 
$S^1\times S^1\overset{id\times p_1}\lo S^1\times \rp^1\overset{f}\lo Y$
extends over the solid torus $S^1\times B^2$.
 \end{proposition}
 \proof Consider the induced map $f':\rp^1\lo Map(S^1,Y)$
to the mapping space defined by $f'(x)(z)=f(z,x)$ for $x\in\rp^1$ and $z\in S^1$.
$f$ extends over $S^1\times \rp^2$ if and only if $f'$ extends
over $\rp^2$. Notice that $f'$ extends over $\rp^2$ if and only if
$S^1\overset{p_1}\lo \rp^1 \overset{f'}\lo Map(S^1,Y)$
extends over $B^2$ which is the same as to say that
$S^1\times S^1\overset{id\times p_1}\lo S^1\times \rp^1\overset{f}\lo Y$
extends over the solid torus $S^1\times B^2$.
\endproof

\begin{proposition}
 \label  {ExtendOverCircleOfRP2} Suppose $(a,b)\in S^1\times \rp^1$.
If $f:S^1\times \rp^1\lo \rp^2$ is a map such that $f|\{a\}\times \rp^1$
is null-homotopic and $f|S^1\times \{b\}$ is not null-homotopic, then $f$
extends over $S^1\times \rp^2$.
 \end{proposition}
 \proof Assume $f|\{a\}\times \rp^1$ is constant.
In view of \ref{ExtendOverTorus} we need to show
that the composition 
$S^1\times S^1\overset{id\times p_1}\lo S^1\times \rp^1\overset{f}\lo \rp^2$
extends over the solid torus $S^1\times B^2$. Let $D$ be a disk with boundary equal to $\rp^1$.
Pick $e:I\lo S^1$ identifying $0$ and $1$ with $a\in S^1$.
The homotopy $f\circ (e\times id):I\times \rp^1\lo \rp^2$
has a lift $H:I\times \rp^1\lo S^2$ such that $\{0\}\times \rp^1$
and $\{1\}\times \rp^1$ are each mapped to a point and those points
are antipodal as $f|S^1\times \{b\}$ is not null-homotopic.
Therefore $H$ can be extended to $G:\partial(I\times D)\lo S^2$
so that $G|\{0\}\times D$ and $G|\{1\}\times D$ are constant. Fix 
the orientation of $\partial(I\times D)$ and let $c$ be the degree of $G$.
Define $F$ on $I\times S^1$ as the composition $G\circ (id\times p_1)$
and use the orientation on $\partial (I\times B^2)$ induced by that on $\partial(I\times D)$.
Define $F$ on $D_1=\{1\}\times B^2$ as a map with the same value on $\partial D_1$
as $G(\partial (\{1\}\times D))$ so that the induced map from $D_1/(\partial D_1)\lo S^2$
is of degree $-c$ (the orientation on $D_1/(\partial D_1)$ is induced
by the orientation of $\partial(I\times B^2)$). The new map is called $F$. 
Define $F$ on $D_0=\{0\}\times B^2$ as the map with the same value on $\partial D_0$
as $G(\partial(\{0\}\times D))$  so that 
$F(0,x)=-F(1,x)$ for all $x\in B^2$.
The cumulative map $F:\partial(I\times B^2)\lo S^2$ is of degree $0$, so it extends
to $F':I\times B^2\lo S^2$. Notice that $J=p_2\circ F':I\times B^2\lo \rp^2$
has the property that $J(0,x)=J(1,x)$ for all $x\in B^2$.
Therefore it induces an extension $S^1\times B^2\lo\rp^2$
of the composition 
$S^1\times S^1\overset{id\times p_1}\lo S^1\times \rp^1\overset{f}\lo \rp^2$.
\qed

\sh{The first modification $M_1$ of $\rp^3$}
 
 Let $B^3 \subset \re^3$ be the unit ball and
 let $D$ be  the  $2$-dimensional disk  of radius $1/3$ lying in
 the $yz$-coordinate plane and centered at the point $(0,1/2,0)$.
 Denote by $L$ the solid torus obtained by rotating $D$
 about the $z$-axis. We consider $L$ with the structure of
cartesian product $L=S^1 \times D$ such that the rotations of $L$ about the $z$-axis
correspond to the rotations of $S^1$. Think of $S^1$ as the circle
 $x^2+y^2=1/4, z=0$  (the circle traced by the center of $D)$. Since $L$ is untouched under the quotient map $q_3:B^3\lo \rp^3$,
we may assume $L\subset \rp^3$.
 The {\it first modification} $M_1$ of $\rp^3$ is obtained
by removing the interior of $L$ from $\rp^3$ and
attaching $S^1\times \rp^2$ via the map $S^1\times \rp^1\to \partial L$,
where $\rp^1$ is identified with $\partial D$.
Notice that $\rp^2=q_3(\partial B^3)\subset M_1$.

  \begin{proposition}
 \label  {p3} There is a retraction $r :M_1 \lo \rp^2$
of the first modification $M_1$ of $\rp^3$ to the projective plane.
 \end{proposition}
 \proof
 We use the  notation
 that we introduced above defining the first modification of $\rp^3$.
   Let $I$ be the interval of the
 points of $B^3$ lying on the $z$-axis and let $M=\partial B^3
 \cup I \subset B^3$.  Denote $K =B^3\sm(L \sm \partial L)$.
Consider the group $\Gamma$ of rotations of $\re^3$ around the $z$-axis.
 Note that $L$, $M$ and $K$ are invariant under  rotations in $\Gamma$
 and every such rotation  induces the corresponding
 homeomorphism of $\rp^2$ which will be called the corresponding
 rotation of $\rp^2$.
 Recall that  $L$ is represented
 as the product $L=S^1\times D$  in such a way
 that the rotations of $L$ are induced by the rotations of  $S^1$.
 Let $\alpha : K \lo M$ be a retraction which commutes with
 the rotations in $\Gamma$ (this means that for every   rotation $\rho\in\Gamma$   of $\re^3$
 and $x\in K$, $\alpha(\rho(x))=\rho(\alpha(x))$).
 Let $\beta : M \lo \rp^2$ be the extension of
 $q_3$ restricted to $\partial B^3$ sending the interval $I$ to
 the point $q_3(\partial I)$.
 Then $\beta$ also commutes with the rotations in $\Gamma$ (this means that
 for every $x \in M$,  every rotation $\rho\in\Gamma$ of $\re^3$
 and the corresponding rotation $\rho'$ of $\rp^2$,
 $\beta(\rho(x))=\rho'(\beta(x))$).
 Denote  $\gamma= \beta \circ\alpha : K
 \lo \rp^2$. It is easy to see
 that  $\gamma$ commutes with the rotations.
 Consider $\partial D$ as a subspace  $\partial D = \rp^1
 \subset \rp^2$ of a projective plane, and
 consider $\partial L=S^1 \times \partial  D$ as the subset of
  $T =S^1 \times \rp^2$ induced by the inclusion $\partial D \subset
  \rp^2$.
Fix $a$ in $S^1$. By Proposition \ref{p0}
               extend $\gamma$ restricted to $\{a\} \times \partial D$
               to a map $\mu : \{a\} \times RP^2 \lo RP^2$ and
               by the rotations of $S^1$ and $\rp^2$
 extend the map $\mu$   to the map $\kappa : T=S^1 \times \rp^2 \lo \rp^2$.
 Note that $\kappa$ agrees with $\gamma$ on $\partial L=S^1 \times \partial D $
 and therefore the maps $\gamma$ and $ \kappa$ define
 the map $\nu :M_1 \lo \rp^2$
 from the first modification of $\rp^3$.
 This map $\nu$ induces  the  required retraction $r : M_1 \lo \rp^2$.
\qed

\sh{The second modification $M_2$ of $\rp^3$}
 
Let $B^3$  be the unit ball in $\re^3$ and
 let $L\subset B^3$ be the subset of  $B^3$ consisting of
 the points lying in the cylinder   $x^2 + y^2 \leq 1/4$.
Notice that $R=q_3(L)$ is a solid torus in $\rp^3$
as the map $B^2\lo B^2$ sending $z$ to $-z$ is isotopic to
the identity. Set $D=R\cap \rp^2$.
Represent $R$ as $S^1\times D$ such that $\{a\}\times D$ 
          is identified with $D$.
The {\it second modification} $M_2$ of $\rp^3$ is obtained
by removing the interior of $R$ and attaching $S^1\times \rp^2\cup \{a\}\times \rp^3$
via the inclusion $S^1\times \partial D\lo S^1\times \rp^2\cup \{a\}\times \rp^3$,
where $\partial D$ is identified with $\rp^1$.

 \begin{proposition}
 \label  {p4} Let $M_2$ be the second modification of $\rp^3$.
 The inclusion $i :\rp^2 \cap M_2  \hookrightarrow
 \rp^2$ extends to a map  $g : M_2 \lo \rp^2$.
 \end{proposition}
 \proof
  We use the  notation
 that we introduced above defining the notion of the second
modification of $\rp^3$.
 Denote $H=\rp^3 \sm (R \sm \partial R)$. Since the center of $B^3$
 does not belong to $q_3^{-1}(H)$, the radial projection sends
 $q_3^{-1}(H)$ to $\partial B^3$ and hence the radial projection
 induces the corresponding map $\alpha : H \lo \rp^2$ which extends
 the map $i$. Recall that $R$ is represented as $S^1\times D$.
  Then $\partial R = S^1 \times \partial D$.
 Fix $a \in S^1$ and $b \in \partial D$.
Notice that $\alpha | \{a\}\times \partial D$ is null-homotopic
and $\alpha | S^1\times \{b\}$ is not null-homotopic.
By \ref{ExtendOverCircleOfRP2} one can extend $\alpha|\partial R$ over $S^1\times \rp^2$.
Any such extension is null-homotopic when restricted to $\{a\}\times\rp^1$, 
so it can be extended over $\{a\}\times\rp^3$ (see \ref{p1}).
Pasting all those extensions gives the desired map $g:M_2\lo \rp^2$.
\endproof
  \end{section}

  \begin{section}{Proof of Theorem \ref{t1}}

\begin{lemma}
 \label  {MapsToCircleOfRP2PlusRP3}
Suppose $X$ is a compactum of dimension at most three and mod 2 dimension
$\dim_{\z_2}X$ of $X$ equals $1$.
A map $f:A\lo S^1\times\rp^2\cup \{a\}\times \rp^3$
extends over $X$ if and only if $\pi\circ f$ extends over $X$,
where $\pi:S^1\times\rp^2\cup \{a\}\times \rp^3\lo S^1$ is the projection
onto the first coordinate.
 \end{lemma}
 \proof Only one direction is of interest.
Pick an extension $g:X\to S^1$ of $\pi\circ f$.
Let $\pi_2:S^1\times \rp^3\lo \rp^3$ be the projection.
Since $\edim X\leq \rp^\infty$ implies $\edim X\leq \rp^3$
(see \ref{skeleta}), the composition $\pi_2\circ f$
extends over $X$ to $h:X\lo \rp^3$.
The diagonal $G:X\lo S^1\times \rp^3$ of $g$ and $h$
can be pushed rel.$A$ to the $3$-skeleton of $S^1\times \rp^3$
which is exactly $S^1\times\rp^2\cup \{a\}\times \rp^3$
under standard CW structures of $S^1$ and $\rp^3$.
The resulting map $X\lo S^1\times\rp^2\cup \{a\}\times \rp^3$
is an extension of $f$.
 \endproof
\begin{corollary}
 \label  {MapsToCircleThenRP2}
Suppose $X$ is a compactum of dimension at most three and $A$ is a closed subset of $X$. If mod 2 dimension
$\dim_{\z_2}X$ of $X$ equals $1$, then any map $f:A\lo \rp^1$ followed by the inclusion
$\rp^1\lo \rp^2$
extends over $X$.
 \end{corollary}
 \proof Let $i:\rp^1\lo \rp^3$ be the inclusion.
Extend $i\circ f:A\lo \rp^3$ to $G:X\lo \rp^3$.
Let $R$ be the solid torus as in the definition of the second
modification of $\rp^3$. Put $Y=G^{-1}(R)$
and $B=G^{-1}(\partial R)$. The map $g:B\to \partial R=S^1\times \rp^1$
induced by $G$ extends to $H:Y\lo S^1\times \rp^2\cup \{a\}\times \rp^3$ in view of 
 \ref{MapsToCircleOfRP2PlusRP3}. Pasting $G|(X\setminus G^{-1}(int(R))$
and $H$ results in an extension $F:X\to M_2$ of $f$.
By \ref{p4} the inclusion $\rp^2 \cap M_2  \hookrightarrow
 \rp^2$ extends to a map  $g : M_2 \lo \rp^2$.
Notice that $g\circ F$ is an extension of $i\circ f$.
 \endproof\begin{corollary}
 \label  {MapsToRP2ThenM1}
Suppose $X$ is a compactum of dimension at most three. If mod 2 dimension
$\dim_{\z_2}X$ of $X$ equals $1$, then any map $f:A\lo \rp^2$ followed by the inclusion
$\rp^2\lo M_1$ from the projective plane to the first modification of $\rp^3$
extends over $X$.
 \end{corollary}
 \proof Let $i:\rp^2\lo \rp^3$ be the inclusion.
Extend $i\circ f:A\lo \rp^3$ to $G:X\lo \rp^3$.
Let $L$ be the solid torus as in the definition of the first
modification of $\rp^3$. Put $Y=G^{-1}(L)$
and $B=G^{-1}(\partial L)$. The map $g:B\to \partial L=S^1\times \rp^1$
induced by $G$ extends to $H:Y\lo S^1\times \rp^2$ in view of 
 \ref{MapsToCircleThenRP2}. Pasting $G|(X\setminus G^{-1}(int(L)))$
and $H$ results in an extension $F:X\to M_1$ of $f$.
 \endproof
Since $\rp^2$ is a retract of $M_1$ (see \ref{p3}),
Corollary \ref{MapsToRP2ThenM1} does indeed imply Theorem \ref{t1}.

\rk{Acknowledgement}This research was supported
by Grant No.\ 2004047  from the United States-Israel Binational Science
Foundation (BSF),  Jerusalem, Israel.

\end{section}

\Addresses\recd

\end{document}